\documentclass{amsart}

\calclayout
\usepackage{amssymb, amsmath}
\usepackage{mathrsfs}
\usepackage{mathtools}
\usepackage{amscd}
\usepackage{verbatim}

\usepackage{wasysym}

\usepackage[utf8]{inputenc}
\usepackage{enumerate}
\usepackage{url}

\usepackage[colorinlistoftodos,prependcaption,textsize=tiny]{todonotes}

\usepackage[colorlinks,linkcolor={blue},citecolor={blue},urlcolor={red},]{hyperref}

\allowdisplaybreaks

\theoremstyle{plain}
\newtheorem{theorem}{Theorem}
\theoremstyle{remark}
\newtheorem{remark}[theorem]{Remark}

\theoremstyle{plain}

\newtheorem{lemma}[theorem]{Lemma}

\numberwithin{theorem}{section}

\numberwithin{equation}{section}

\def\N{{\mathbb N}}

\def\R{{\mathbb R}}

\newcommand{\E}{{\mathbb E}}
\renewcommand{\P}{{\mathbb P}}

\newcommand{\mbf}[1]{\mathbf{#1}}

\begin{document}

\title[Branching Process Representation for PDEs]{Branching process representation for nonlinear first-order conservation PDEs in 1D}

\author{Jochem Hoogendijk}\email{j.p.c.hoogendijk@uu.nl}
\address{Mathematical institute, 
    Utrecht University, 3508 TA Utrecht, The
    Netherlands} 

\author{Ivan Kryven}\email{i.v.kryven@uu.nl}
\address{Mathematical institute,
    Utrecht University, 3508 TA Utrecht, The
    Netherlands} 


\begin{abstract}
    We show that a large class of 1D first-order conservation PDEs can be probabilistically represented using multi-type branching processes. The representation holds when the initial conditions are linear combinations of negative exponentials. We also show that in some cases, the time of gradient blow up can be identified by studying criticality conditions of the corresponding branching processes.
\end{abstract}

\maketitle

\section{Introduction}
\label{sec:intro}
This paper introduces a Feynman-Kac-type formula for the class of first-order conservation equations
\begin{equation}\label{eq:first_order_cons_eq_main}
  u_t + (F(u))_x = 0, \qquad u(0, x) = \sum\limits_{k=0}^m p_k e^{-kx}, \qquad x \in \R,
\end{equation}
with detailed assumptions given in Section \ref{sec:main_results}. Roughly speaking, we show that $u(t, x) = \E[g(x, \mbf{T})],$ where $\mbf{T}$ is the total progeny vector of a multi-type branching process and $g$ is some function we will specify later in the paper. The branching process incorporates both the initial condition $u(0, x)$ and the nonlinearity $F$, and incorporates time $t$ as a parameter. We call such a solution the \emph{branching process representation} of the solution of \eqref{eq:first_order_cons_eq_main}. 

A particularly appealing feature of the branching process representation is that it can capture infinite gradients that may arise in the solution of the PDE\eqref{eq:first_order_cons_eq_main} under certain initial conditions and nonlinearities. More precisely, it is the criticality of the branching process that determines when and where the solution of the PDE \eqref{eq:first_order_cons_eq_main} develops infinite gradients.

Analogous to how the heat equation plays a central role for the original Feynman-Kac formula, our technique is inspired by the inviscid Burgers' equation $u_t + u u_x = 0$, which is arguably the simplest first-order conservation PDE. Burgers' equation, when supplied with very specific initial conditions, becomes the  generating function transform of Smoluchowski's multiplicative coalescence equation \cite{menon2004,deaconu2000}. The latter has been shown to be closely related to the connected components statistics in random graphs \cite{aldous_1999,normand2011,bertoin2012}, which can also be studied with branching processes through the so-called exploration principle \cite{durrett_2007,bertoin2010}. This suggests that one may read off the solution of Burgers' equation from the total progeny distribution of a branching process with a certain offspring distribution. To expand this observation into a technique, one needs to account for general initial conditions and general nonlinearity terms. In this paper, we show that switching to more general multi-type branching processes resolves both issues. To prove this statement we take a simple analytic route that does not depend on the analogy with random graphs. We use distributional self-similarity of the multi-type branching process together with Chernoff's bound to show that the proposed solution is well-defined and satisfies the PDE \eqref{eq:first_order_cons_eq_main}.

Apart from the Feynman-Kac formula (also known as the Kolmogorov forward/backward equation \cite{revuz_yor_1999, durrett_2010}), other probabilistic interpretations of solutions of PDEs also appear in the literature. For instance, the asymmetric simple exclusion process (ASEP) is an interacting particle system that has the inviscid Burgers' equation as its hydrodynamic limit \cite{ferrari_1992,rezakhanlou_1995}. There is also the application of branching Brownian motion to study the Fisher-Kolmogorov-Petrovskii-Piskunov (FKPP) equation and related semilinear parabolic generalizations \cite{mckean_1975, an_2023}.

The rest of the paper is structured as follows. In Section \ref{sec:prelim} we introduce some notation and definitions as well as the distributional self-similarity lemma that will be used throughout the rest of the paper. Section \ref{sec:main_results} contains the two main results in this paper, as well as some discussion of them. Section \ref{sec:proofs} contains the proofs of the main results.

\section{Preliminaries}
\label{sec:prelim}
\subsection{Notation}
Given $m \in \N$, we define $[m] := \{1, \ldots, n\}$ and $[m]_0 := \{0\} \cup [m]$.

\subsection{Multi-type branching processes}
We define the multi-type branching process as follows. For fixed type $k \in [m]_0$, let $\mathbf{X}_{k,(n,i)}$ with $(n, i) \in \N_0 \times \N$ be a collection of independent and identically distributed copies of a non-negative discrete random vector $\mathbf{X}_{k} = (X_{k, 0}, \ldots, X_{k, m})$. The \emph{multi-type branching process started from type $k$} is the recurrence equation 
\begin{equation*}
    \mathbf{Z}_{n+1} = \sum\limits_{k=0}^m \sum\limits_{i=1}^{Z_{n, k}}\mathbf{X}_{k,(n, i)},
\end{equation*}
initialized with the standard unit vector, $\mathbf{Z}_0 = \mathbf{e}_k$. Then, the random vector $\mathbf{T}^{(k)} = \sum_{n=0}^\infty \mathbf{Z}_n$ is called the \emph{total progeny of the multi-type branching process started from type $k$}. Note that given $j \in [m]_0$, $T_j^{(k)}$ is the number of individuals of type $j$ in a multi-type branching process started from type $k$. 

If the starting type $k$ is chosen at random with probability $p_k$, i.e. $\P(\mathbf{Z}_0 = \mathbf{e}_k) = p_k$, we define the \emph{overall total progeny} $\mathbf{T}$ as $\mbf{T} = \sum_{n=0}^\infty \mathbf{Z}_n$. For more information on branching processes see \cite{harris_1963, athreya_1972}.

\subsection{Probability generating functions}
Given a random variable $\mathbf{X}$ taking values in $\N_0^m$ its generating function is defined as $$G_{\mathbf{X}}(s_1, \ldots, s_m) = \sum\limits_{\mathbf{n} \in \N_0^m} \P(\mathbf{X} = \mathbf{n}) s_1^{n_1} \cdot \ldots \cdot s_m^{n_m}.$$ 
One particularly useful property in the case of multi-type branching processes is the following lemma:
\begin{lemma}[\cite{good_1955}]\label{lem:total_progeny_implicit}
    Given a multi-type branching process with independent and identically distributed offspring $\mathbf{X}_k$ for each type $k \in [m]_0$ with $m \in \N_0$, the following formal implicit system of equations holds:
    \begin{equation}\label{eq:total_progeny_multi_type}
        \begin{split}
                G_{\mbf{T}^{(0)}} &= s_0 G_{\mbf{X}_0}(G_{\mbf{T}^{(0)}}, \ldots, G_{\mbf{T}^{(m)}}),\\
                &\stackrel{\vdots}{\phantom{=}}\\
                G_{\mbf{T}^{(m)}} &= s_m G_{\mbf{X}_m}(G_{\mbf{T}^{(0)}}, \ldots, G_{\mbf{T}^{(m)}}),
        \end{split}
    \end{equation}
    where $\mathbf{s} \in \mathbb{K}^{m+1}$ is a vector of formal variables and $\mathbb{K}$ is a field.
\end{lemma}

\section{Main results}
\label{sec:main_results}
Throughout the rest of this paper, we will consider the following. Let $N, m \in \N_0$ and consider the PDE
\begin{equation}\label{eq:first_order_cons_pde}
    u_t + F'(u) u_x = 0, \qquad u(0, x) = \sum\limits_{k=0}^m p_k e^{-kx}, \qquad x \in \R
\end{equation}
where $F$ is of the form $$F(s) = \sum\limits_{n=0}^{N+1} f_n s^n$$ with $f_n, p_k \in \R$ for all $n \in [N]_0$ and $k \in [m]_0$ respectively. 

We start with a special case of our main theorem.
\begin{theorem}\label{thm:main_result_2}
    Consider the PDE \eqref{eq:first_order_cons_pde} with $(p_k)_{k=0}^m$ being a probability distribution and $(f_n)_{n=0}^{N+1}$ such that $\sum_{n=0}^{N+1} n f_n = 0$, $f_1 \geq -1$ and $f_n \geq 0$ for all $n \geq 2$. Let $[m]_0$ be the type space of a branching process with offspring distribution determined by
    \begin{equation}
        G_{\mbf{X}_k}(\mbf{s}) = \exp\left(tkF'\left(\sum\limits_{l=0}^m p_l s_l\right)\right).
    \end{equation}
    Then, the solution of the PDE \eqref{eq:first_order_cons_pde} is given by
    \begin{equation}\label{eq:solution_prob_setting}
        u(t, x) = G_{\mbf{T}}(e^{-x}, \ldots, e^{-mx}), \qquad \text{ on } (0, t_c) \times [0, \infty),
    \end{equation}
    with $t_c$ given by
    \begin{equation}\label{eq:critical_time}
        t_c = \left(\sum\limits_{k=0}^m k p_k\right)^{-1} \left(\sum\limits_{n=2}^{N+1} n(n-1) f_n\right)^{-1}
    \end{equation}
    at which a gradient blow-up occurs at $x = 0$.
\end{theorem}
\begin{remark}
    The solution \eqref{eq:solution_prob_setting} can be seen as a Feynman-Kac type formula because it can be written as
    \begin{equation}
        u(t, x) = \E\left[\exp\left(-x\left(\sum\limits_{k=0}^m k T_k\right)\right)\right].
    \end{equation}
\end{remark}
\begin{remark}
    The assumption on the coefficients $(f_n)$ is equivalent to $F'(s) + 1$ being a probability generating function.
\end{remark}
\begin{remark}
    By the method of characteristics we have uniqueness, so it indeed makes sense to speak of \emph{the} solution.
\end{remark}
\begin{remark}
    Note that we pose the PDE on $\R$ and consider a representation of the solution on $[0, \infty)$, meaning that we do not have to specify a boundary condition at $x = 0$. Alternatively, one can pose the PDE on $[0, \infty)$ and supply it with the corresonding boundary condition induced by $u(0, x)$ for $x < 0$.
\end{remark}
Theorem \ref{thm:main_result_2} says that when the initial condition and nonlinearity have probabilistic structure, the gradient blow-up occurs at the time when the multi-type branching process becomes critical, i.e. $\E[|\mbf{T}|] = \infty$. In this case, the branching process provides a direct way to deduce when a classical solution of a PDE fails to exist.

The branching process representation in Theorem \ref{thm:main_result_2} can be generalized by dropping the assumptions on the coefficients $(p_k)_{k=0}^m$ and $(f_n)_{n=0}^{N+1}$. In this case, the solution of the PDE \eqref{eq:first_order_cons_pde} can still be expressed in terms of a multi-type branching process, but more types have to be introduced, as well as an auxiliary probability distribution $(q_k)_{k=0}^{Nm}$. Before introducing the theorem itself, we introduce some notation combining the coefficients of the nonlinearity and the initial condition. First, define $\Phi(s) = F(s) + s$, so that $\Phi = \sum_{n=0}^{N+1} \phi_n s^n$ with $\phi_n = f_n$ for $n \neq 1$ and $\phi_1 = f_1 + 1$. Then, define the coefficients $(H_k)_{k=0}^{Nm}$ as those resulting from the polynomial expansion
\begin{equation}\label{eq:h_k_expansion}
    \sum\limits_{n=0}^{N+1} n \phi_n \left(\sum\limits_{k=0}^m p_k s^k\right)^{n-1} =: \sum\limits_{k=0}^{Nm} H_k s^k.
\end{equation}
The following theorem holds.
\begin{theorem}\label{thm:main_result_1}
    Consider a multi-type branching process on the type space $[Nm]_0$ with offspring distribution determined by 
    \begin{equation}\label{eq:offspring_distribution}
        G_{\mbf{X}_k}(\mbf{s}) = \exp\left(tk\left(\sum\limits_{l=0}^{Nm} q_l s_l - 1\right)\right)
    \end{equation} 
    for type $k \in [Nm]_0$, with $(q_l)_{l \in [Nm]_0}$ being an auxiliary probability distribution fully supported on $[Nm]_0$. Furthermore, assume $H_k \neq 0$ for all $k \in [Nm]_0$. Then, there exists $\tau > 0$, such that the solution of the PDE \eqref{eq:first_order_cons_pde} on $(0, \tau) \times [0, \infty)$ is given by
    \begin{equation}\label{eq:pde_solution}
        u(t, x) = \sum\limits_{k=0}^m \frac{p_k q_k}{H_k} G_{\mbf{T}^{(k)}}\left(\frac{H_0}{q_0}, \frac{H_1}{q_1} e^{-x}, \ldots, \frac{H_{Nm}}{q_{Nm}} e^{-Nmx}\right).
    \end{equation}
\end{theorem}
\begin{remark}
    The assumption that $H_k \neq 0$ for all $k \in [Nm]_0$ is only made for expository purposes, but without loss of generality. In case $H_k = 0$ for some $k \in [Nm]_0$, the type space needs to be modified by removing type $k$.
\end{remark}
The above theorem says that the solution of the PDE \eqref{eq:first_order_cons_pde} is given by the weighted expectation value of some functional of a multi-type branching process. That is,
\begin{equation}
    u(t, x) =  \mathbb{E}\left[\sum\limits_{k=0}^m \frac{p_k q_k}{H_k} \left(\frac{H_0}{q_0}\right)^{T_0^{(k)}}\cdot \left(\frac{H_1}{q_1} e^{-x}\right)^{T_1^{(k)}} \cdot \ldots \cdot \left(\frac{H_{Nm}}{q_{Nm}} e^{-Nm x}\right)^{T_{Nm}^{(k)}}\right]
\end{equation}
Therefore, the solution $u(t, x)$ is determined by the probabilities $\P(\mbf{T}^{(k)} = \mbf{n})$ associated with the total progenies of the multi-type branching process. Furthermore, note that given a type $k \in [Nm]_0$, the offspring distribution of $X_{k, j} \sim \mathrm{Poi}(tk q_j)$ for $j \in [Nm]_0$. This observation can be used to show that 
\begin{equation}
    \P(\mbf{T}^{(k)} = \mbf{n}) = e^{-tM(\mbf{n})} \frac{(tM(\mbf{n})q_0)^{n_0}}{n_0!} \cdot \ldots \cdot \frac{(tM(\mbf{n})q_{Nm})^{n_{Nm}}}{n_{Nm}!} \frac{k n_k} {t M(\mbf{n})^2 q_k},
\end{equation}
where $M(\mbf{n}) = \sum_{l=0}^{Nm}ln_l$, by applying the Lagrange inversion formula \cite{good_1960}.

By choosing some specific probability distribution $(q_k)_{k=0}^{Nm}$, we can, in principle, always use the same branching process, but the argument for which we evaluate the generating function changes based on the initial condition and the nonlinearity. As we will see from the proof, the only role of the auxiliary probability distribution is to ensure that we are dealing with \emph{probability} generating functions, which allows us to make use of probabilistic tools.

Finally, note that the theorem only states that the branching process representation of the solution of the PDE \eqref{eq:first_order_cons_pde} is valid at least on $(0, \tau) \times [0, \infty)$, but the PDE itself was posed on $\R_{\geq 0} \times \R$. The reason is that the branching process representation $u(t, x)$ is expressed as a power series in terms of the variable $e^{-x}$. Since the power series we consider have a finite radius of convergence, we cannot guarantee that the branching process representation $u(t, x)$ is well-defined for all $t \geq 0$.

\section{Proofs}
\label{sec:proofs}
This section contains two parts. Both contain the main ideas and proof of the Theorems \ref{thm:main_result_1} and \ref{thm:main_result_2}, respectively. Note that the theorems are proven in opposite order.
\subsection{Proof of Theorem \ref{thm:main_result_1}}
The idea of the proof is as follows. First, we have to show that the candidate solution $u(t, x)$ as well as its derivatives are well-defined. This is done by bounding the power series by a single-type generating function of the total progeny. This way, we can easily obtain explicit information on the coefficients by using Lagrange inversion and Chernoff's inequality to obtain the tail behavior. By choosing $\tau$ small enough, we can ensure well-definedness. The second step is to show that the candidate solution $u(t, x)$ actually solves the PDE \eqref{eq:first_order_cons_pde}. This is done by establishing an implicit equation for $u(t, x)$ based on the branching process properties, and then using implicit differentiation to show that $u(t, x)$ solves the PDE.
\begin{proof}[Proof of Theorem \ref{thm:main_result_1}]
    The proof consists of two steps. First, we have to show that there exists $\tau > 0$ such that $u(t, x)$ and its derivatives, are well-defined on $(0, \tau) \times [0, \infty)$. Notice that well-definedness has to be shown because $u(t, x)$ is expressed as a power series in terms of the variable $e^{-x}$. The second and final step consists of actually showing that the candidate solution $u(t, x)$ satisfies the PDE \eqref{eq:first_order_cons_pde}.

    \textit{Step 1: Showing that $u$, $u_x$ and $u_t$ are well-defined}: Notice that there exist $M > 0$ and $c > 0$ such that
    \begin{equation}
        |u(t, x)| \leq M \sum\limits_{k=0}^{Nm} q_k G_{\mbf{T}^{(k)}}(c, \ldots, c) = M G_{\mbf{T}}(c, \ldots, c).
    \end{equation}
    Define the formal power series $G_{|\mbf{T}|}(s) := G_{\mbf{T}}(s, \ldots, s)$ and $G_X(s) = \sum_{k=0}^{Nm} q_k \exp(tk(s - 1))$, where $|\mbf{T}| := T_0 + \ldots + T_{Nm}$. Then, we see that by using the implicit equation \eqref{eq:total_progeny_multi_type} for $G_{\mbf{T}^{(k)}}$ and the offspring distribution PGF \eqref{eq:offspring_distribution},
    \begin{equation}
        \begin{split}
            G_{|\mbf{T}|}(c) &= \sum\limits_{k=0}^{Nm} q_k G_{\mbf{T}^{(k)}}(c, \ldots, c)\\
            &= c \sum\limits_{k=0}^{Nm} q_k \exp\left(tk\left(\sum\limits_{l=0}^{Nm} q_l G_{\mbf{T}^{(l)}}(c, \ldots, c) - 1\right)\right)\\
            &= c G_X(G_{|\mbf{T}|}(c)).
        \end{split}
    \end{equation}
    Observe therefore that $|\mbf{T}|$ is the total progeny of a single-type Galton-Watson branching process with offspring distribution $X$. Applying the Lagrange inversion formula to $G_{|\mbf{T}|}(c)$, we obtain that
    \begin{equation}\label{eq:gen_cons_eq_lag_inv}
        \P(|\mbf{T}| = n) = \frac{1}{n}\P(X_1 + \ldots + X_n = n - 1),
    \end{equation}
    where $(X_i)$ are i.i.d. copies of the random variable $X$. The right-hand side (RHS) can be estimated using Chernoff's inequality. Letting $\lambda \in \R$ be arbitrary, we obtain the bound 
    \begin{equation}\label{eq:gen_cons_eq_chernoff}
        \begin{split}
            \P(|\mbf{T}| = n) c^n &\leq \frac{1}{n} e^{-\lambda (n-1)} c^n \E[e^{\lambda X}]^n\\
            &= \frac{1}{n} e^{\lambda} \left(c \sum\limits_{k=0}^{Nm} e^{-\lambda} q_k \exp(tk(e^{\lambda} - 1))\right)^n.
        \end{split}
    \end{equation}
    Notice that for any $\lambda \in \R$, the term $\exp(tk(e^{\lambda} - 1))$ can be made as close to 1 as desired by possibly choosing $t > 0$ small. Therefore, there exists $\tau_1 > 0$ and  $\lambda^* \in \R$ such that for all $t \in (0, \tau_1)$, 
    \begin{equation}\label{eq:gen_cons_eq_chernoff_criterion}
        \sum\limits_{k=0}^{Nm} e^{-\lambda^*} q_k \exp(tk(e^{\lambda^*} - 1)) c < 1.
    \end{equation}
    We therefore conclude that for $t \in (0, \tau_1)$,
    \begin{equation}
        G_{|\mbf{T}|}(c) = \sum\limits_{n=1}^\infty \P(|\mbf{T}| = n) c^n < \infty,
    \end{equation}
    which implies that $u(t, x)$ is well-defined on $(0, \tau_1) \times [0, \infty)$.
    
    Next, we show that there exists $\tau_2$ such that $u_x$ is well-defined on $(0, \tau_2) \times [0, \infty)$. By the Leibniz rule, it suffices to show that
    \begin{equation}
        \sum\limits_{k=0}^m \frac{|p_k q_k|}{|H_k|} \mathbb{E}\left[\left(\sum\limits_{l=0}^{Nm} l T_l^{(k)}\right) \left(\frac{|H_0|}{q_0}\right)^{T_0^{(k)}} \cdot \ldots \cdot \left(\frac{|H_{Nm}|}{q_{Nm}}\right)^{T_{Nm}^{(k)}}\right] < \infty.
    \end{equation}
    By the Cauchy-Schwarz inequality, it suffices in turn to show that for each $k \in [m]_0$ that
    $$\E\left[\left(\sum\limits_{l=0}^{Nm} l T_l^{(k)} \right)^2 \right] < \infty \quad \text{ and } \quad \E\left[\left(\frac{|H_0|}{q_0}\right)^{2T_0^{(k)}} \cdot \ldots \cdot \left(\frac{|H_{Nm}|}{q_{Nm}}\right)^{2T_{Nm}^{(k)}}\right] < \infty.$$ By following a similar calculation as before and choosing $\tau_2$ appropriately, it follows that the second quantity is well-defined on $(0, \tau_2) \times [0, \infty)$. The first quantity is finite since in the subcritical regime, the tails of $\mbf{T}^{(k)}$ decrease exponentially as follows from \eqref{eq:gen_cons_eq_chernoff} and \eqref{eq:gen_cons_eq_chernoff_criterion}.

    Finally, we show that $u_t$ is well-defined. By applying the Leibniz rule twice, it suffices to show that $\frac{\partial}{\partial t} G_T(c) < \infty$. The Implicit Function Theorem applied to $G_T(c) = c G_X(G_T(c))$ shows that $G_T(c)$ is differentiable with respect to $t$ on $(0, \tau_3)$ for some $\tau_3 > 0$.

    This finishes the proof that $u(t, x)$ and its derivatives are well-defined on $(0, \tau) \times [0, \infty)$ for some $\tau > 0$, by choosing 
    \begin{equation} 
        \tau_4 := \min\{\tau_1, \tau_2, \tau_3\}.
    \end{equation}

    \textit{Step 2: Showing that \eqref{eq:pde_solution} solves \eqref{eq:first_order_cons_pde}}: We will first show that $u(t, x)$ satisfies a specific implicit equation. This implicit equation will later help us to show that $u(t, x)$ solves the PDE \eqref{eq:first_order_cons_pde}. Henceforth, we shall suppress the arguments of $G_{\mbf{T}^{(k)}}$ by writing 
    \begin{equation}
        G_{\mbf{T}{^{(k)}}} := G_{\mbf{T}^{(k)}}\left(\frac{H_0}{q_0}, \frac{H_1}{q_1} e^{-x}, \ldots, \frac{H_{Nm}}{q_{Nm}} e^{-Nm x}\right).
    \end{equation}

    First, we have
    \begin{equation}\label{eq:first_implicit_relation}
        \begin{split}
            u(t, x) &\stackrel{(i)}{=} \sum\limits_{k=0}^m \frac{p_k}{H_k} q_k G_{\mbf{T}^{(k)}}\left(\frac{H_0}{q_0}, \frac{H_1}{q_1} e^{-x}, \ldots, \frac{H_{Nm}}{q_{Nm}} e^{-Nmx}\right)\\
            &\stackrel{(ii)}{=} \sum\limits_{k=0}^m p_k e^{-k x} G_{\mbf{X}_k}(G_{\mbf{T}^{(0)}}, \ldots, G_{\mbf{T}^{(Nm)}})\\
            &\stackrel{(iii)}{=} \sum\limits_{k=0}^m p_k e^{-k x} \exp\left(tk \left(\sum\limits_{l=0}^{Nm} q_l G_{\mbf{T}^{(l)}} - 1\right)\right)\\
            &\stackrel{(iv)}{=}\sum\limits_{k=0}^m p_k e^{-k x} \exp\left(tk \left(G_{\mbf{T}} - 1\right)\right),
        \end{split}
    \end{equation}
    where $(i)$ follows from the definition \eqref{eq:pde_solution}, $(ii)$ follows by Lemma \ref{lem:total_progeny_implicit}, $(iii)$ follows from the definition of the offspring distribution \eqref{eq:offspring_distribution} and $(iv)$ follows from the definition of $G_{\mbf{T}}$. To further manipulate the above expression, we also deduce an implicit equation for $G_{\mbf{T}}$. By the implicit equation \ref{lem:total_progeny_implicit} for the total progeny of the multi-type branching process and the definition of the coefficients $(H_k)$ \eqref{eq:h_k_expansion}, we have that
    \begin{equation}\label{eq:second_implicit_relation}
        \begin{split}
            G_{\mbf{T}} = \sum\limits_{k=0}^{Nm} q_k G_{\mbf{T}^{(k)}} &= \sum\limits_{k=0}^{Nm} H_k e^{-kx} \exp\left(tk(G_{\mbf{T}} - 1)\right)\\
            &= \sum\limits_{n=0}^N n \phi_n \left(\sum\limits_{k=0}^m p_k e^{-kx} \exp(tk(G_{\mbf{T}} - 1))\right)^{n-1}.
        \end{split}
    \end{equation}
    Combining \eqref{eq:first_implicit_relation} and \eqref{eq:second_implicit_relation} gives
    \begin{equation}\label{eq:implicit_solution}
        \begin{split}
            u(t, x) &= \sum\limits_{k=0}^m p_k e^{-kx} \exp\left(tk (G_{\mbf{T}} - 1)\right)\\
            &= \sum\limits_{k=0}^m p_k e^{-kx} \exp\left(tk \left(\sum\limits_{n=0}^N n \phi_n \left(\sum\limits_{l=0}^{m} p_l e^{-lx} \exp(tl (G_{\mbf{T}^{(l)}} - 1))\right)^{n-1} - 1\right)\right)\\
            &= \sum\limits_{k = 0}^m p_k e^{-kx} \exp\left(tk \left(\phi'(u(t, x))  - 1\right)\right)\\
            &= \sum\limits_{k = 0}^m p_k e^{-kx} \exp\left(tk F'(u(t, x))\right).
        \end{split}
    \end{equation}
    Having established the implicit equation \eqref{eq:implicit_solution} for the candidate solution $u(t, x)$, we apply implicit differentiation to show that $u(t, x)$ solves the PDE \eqref{eq:first_order_cons_pde}. Denoting 
    \begin{equation}
        A(t, x) := 1 - t \sum\limits_{k=0}^m k p_k e^{-kx} \exp\left(tk \left(F'(u(t, x)) - 1\right)\right) F''(u(t, x)),
    \end{equation}
    we see that
    \begin{equation}
        u_t = \frac{1}{A(t, x)}\left(\sum\limits_{k=0}^m k p_k e^{-kx} \exp\left(tk \left(F'(u(t, x)) - 1\right)\right)\right) (F'(u(t, x)) - 1)
    \end{equation}
    and
    \begin{equation}
        u_x = -\frac{1}{A(t, x)}\left(\sum\limits_{k=0}^m k p_k e^{-kx} \exp\left(tk \left(F'(u(t, x)) - 1\right)\right)\right).
    \end{equation}
    Since we can always choose $\tau < \tau_4$ small enough such that $A(t, x) > 0$ for all $(t, x) \in (0, \tau) \times [0, \infty)$, we have shown that $u(t, x)$ satisfies $$u_t + F'(u) u _x = 0.$$ Hence we have shown that \eqref{eq:pde_solution} solves the PDE \eqref{eq:first_order_cons_pde}.
\end{proof}

\subsection{Proof of Theorem \ref{thm:main_result_2}}
\begin{proof}
    First, observe that the branching process representation $u(t, x)$ is well-defined on $[0, \infty)$ for all $t \geq 0$. Moreover, $u(t, x)$ is decreasing in $x$. Next up are the derivatives of $u(t, x)$. Notice that 
    \begin{equation}
        \frac{\partial}{\partial x} u = \E\left[\left(\sum\limits_{k=0}^m k T_k\right) e^{-x \left(\sum\limits_{k=0}^m k T_k\right)}\right],
    \end{equation}
    meaning that $\frac{\partial}{\partial x} u(t, x)$ is well-defined on $[0, \infty)$ as long as the branching process is subcritical. For $\frac{\partial}{\partial t} u(t, x)$, applying the Leibniz rule twice shows that $\frac{\partial}{\partial t} u(t, x)$ is well-defined on $[0, \infty)$ as long as the branching process is subcritical.
    
    We also see that $u(t, x)$ satisfies the implicit equation
    \begin{equation}
        u(t, x) = \sum\limits_{k=0}^m p_k e^{-kx} \exp\left(tkF'(u(t, x))\right).
    \end{equation}
    Applying implicit differentiation as in the previous section, we see that $u(t, x)$ solves the PDE \eqref{eq:first_order_cons_pde} as long as 
    \begin{equation}
        A(t, x) := 1- t\sum\limits_{k=0}^m k p_k e^{-kx} \exp\left(tkF'(u(t, x))\right) F''(u(t, x)) \neq 0.
    \end{equation}
    Since $u(t, x)$ is decreasing in $x$, we see that for all $t$ such that $A(t, x) \neq 0$, we have that
    \begin{equation}
        A(t, x) \geq A(t, 0).
    \end{equation}
    Therefore, as soon as $A(t, 0) = 0$, a gradient blow-up occurs at $x = 0$. Solving $A(t, 0) = 0$ gives the critical time $t_c$ as in \eqref{eq:critical_time}, which is also the time at which the branching process becomes critical.
\end{proof}

\section*{Acknowledgements}
The authors greatly acknowledge discussions with Rik Versendaal, which helped the project in its early phase. The authors thank Niek Mooij and Mike de Vries for numerous stimulating discussions.
This publication is part of the project ``Random graph representation of nonlinear evolution problems"  of the research programme Mathematics Cluster/NDNS+ which is partly financed by the Dutch Research Council (NWO). IK gratefully acknowledges support form Netherlands Research Organisation (NWO), research program VIDI, project number VI.Vidi.213.108.


\bibliographystyle{alpha}
\bibliography{references.bib}

\end{document}